# A STRONG LAW OF LARGE NUMBERS FOR CAPACITIES[1]

By Fabio Maccheroni and Massimo Marinacci

*Università Bocconi and Università di Torino*

We consider a totally monotone capacity on a Polish space and a sequence of bounded p.i.i.d. random variables. We show that, on a full set, any cluster point of empirical averages lies between the lower and the upper Choquet integrals of the random variables, provided either the random variables or the capacity are continuous.

**1. Introduction.** In this paper we prove a strong law of large numbers for totally monotone capacities. Specifically, given a totally monotone capacity $\nu$ defined on the Borel $\sigma$-algebra $\mathcal{B}$ of a Polish space $\Omega$, and a sequence $\{X_n\}_{n \geq 1}$ of bounded, pairwise independent and identically distributed random variables, we show that

$$\nu\bigg(\bigg\{\omega \in \Omega : \int X_1 \, d\nu \leq \liminf_n \frac{\sum_{j=1}^n X_j(\omega)}{n} \\ \leq \limsup_n \frac{\sum_{j=1}^n X_j(\omega)}{n} \leq -\int -X_1 \, d\nu\bigg\}\bigg) = 1,$$

provided the $X_n$s are continuous or simple, or $\nu$ is continuous. In this way we extend earlier results of Marinacci [13].

Under different names, totally monotone capacities have been widely studied in both pure and applied mathematics. They have been introduced by Choquet [4] motivated by some problems in potential theory, and in his wake many works have studied them in both potential theory and probability theory (see, e.g., [6] and [8]).

In mathematical statistics and in mathematical economics, totally monotone capacities have been used to represent subjective prior beliefs when the information on which such beliefs are based is not good enough to represent them by a standard additive probability (see, e.g., [7, 11, 15, 20] and [19]).

Received February 2004; revised June 2004.

[1]Supported by the Ministero dell'Istruzione, dell'Università e della Ricerca.

Our result shows that even in a nonadditive setting, the limit behavior of empirical averages has some noteworthy properties. In particular, we show that eventually empirical averages lie, with probability one, between the lower and upper Choquet integrals associated with the given capacity. This extends to the nonadditive setting the classic Kolmogorov limit law, to which our result reduces when $\nu$ is additive since in this case lower and upper Choquet integrals coincide.

In a subjective probability perspective, our result says that, while in the additive case a Bayesian decision maker believes that the empirical averages of a sequence of p.i.i.d. random variables tend to a given number, here he just believes that the limit behavior of the empirical averages is confined in a given interval. This reflects a possible lack of confidence in his probability assessments. This interpretation of nonadditive limit laws and its relevance in mathematical economics has been recently discussed at length in [9], to which we refer the interested reader for details and references.

**2. Preliminaries.** Let $\Omega$ be a Polish space and $\mathcal{B}$ its Borel $\sigma$-algebra. A *random variable* (r.v.) is a (Borel) measurable function $X : \Omega \to \mathbb{R}$. A *totally monotone capacity* on $\mathcal{B}$ is a set function $\nu : \mathcal{B} \to [0,1]$ such that:

(c.1) $\nu(\varnothing) = 0$ and $\nu(\Omega) = 1$,
(c.2) $\nu(A) \leq \nu(B)$ for all Borel sets $A \subseteq B$,
(c.3) $\nu(B_n) \downarrow \nu(B)$ for all sequences of Borel sets $B_n \downarrow B$,
(c.4) $\nu(G_n) \uparrow \nu(G)$ for all sequences of open sets $G_n \uparrow G$,
(c.5) $\nu(\bigcup_{j=1}^n B_j) \geq \sum_{\varnothing \neq J \subseteq \{1,\ldots,n\}} (-1)^{|J|+1} \nu(\bigcap_{j \in J} B_j)$ for every collection $B_1, \ldots, B_n$ of Borel sets.

A set function $\nu : \mathcal{B} \to [0,1]$ such that:

(c.6) $\nu(B_n) \uparrow \nu(\Omega)$ for all sequences of Borel sets $B_n \uparrow \Omega$,

is called *continuous*. A continuous set function $\nu : \mathcal{B} \to [0,1]$ is a totally monotone capacity if and only if (c.1), (c.2) and (c.5) hold (see [18] and [14], Theorem 10).

Let $\nu$ be a totally monotone capacity on $\mathcal{B}$. As in the additive case, we say that the elements of a sequence $\{X_n\}_{n \geq 1}$ of r.v.s are *pairwise independent* with respect to $\nu$ if, for each $n, m \geq 1$ and for all open subsets $G_n, G_m$ of $\mathbb{R}$,

$$\nu(\{X_n \in G_n, X_m \in G_m\}) = \nu(\{X_n \in G_n\})\nu(\{X_m \in G_m\});$$

we say that they are *identically distributed* if, for each $n, m \geq 1$ and each open subset $G$ of $\mathbb{R}$,

$$\nu(X_n \in G) = \nu(X_m \in G).$$

The *Choquet integral* of a bounded r.v. $X$ with respect to a totally monotone capacity $\nu$ is defined by

$$\int X \, d\nu \equiv \int_0^{+\infty} \nu(\{X > t\}) \, dt + \int_{-\infty}^0 [\nu(\{X > t\}) - 1] \, dt.$$

The integrals on the right-hand side are Riemann integrals and they are well defined since $\nu(\{X > t\})$ is a monotone function in $t$. The Choquet integral is positively homogeneous, monotone and translation invariant [i.e., $\int (X+c)\,d\nu = \int X\,d\nu + c$ if $c$ is constant]. It reduces to the standard integral when $\nu$ is an additive probability measure.

In general, $\int X\,d\nu \leq -\int -X\,d\nu$. Equality holds for all r.v.s if and only if $\nu$ is additive. The integrals $\int X\,d\nu$ and $-\int -X\,d\nu$ are sometimes called *lower* and *upper Choquet integrals*, respectively.

**3. The law of large numbers.** We can now state our main result.

THEOREM 1. *Let $\nu$ be a totally monotone capacity on $\mathcal{B}$, and $\{X_n\}_{n\geq 1}$ a sequence of bounded, pairwise independent and identically distributed random variables. Then*

$$\nu\left(\left\{\omega \in \Omega : \int X_1\,d\nu \leq \liminf_n \frac{\sum_{j=1}^n X_j(\omega)}{n}\right.\right.$$
$$\left.\left.\leq \limsup_n \frac{\sum_{j=1}^n X_j(\omega)}{n} \leq -\int -X_1\,d\nu\right\}\right) = 1,$$

*provided at least one the following two conditions holds*:

(i) *$\nu$ is continuous;*
(ii) *the random variables $X_n$ are either continuous or simple.*

A few remarks are in order. First, as $\nu(B) = 1$ and $A \subseteq B^c$ imply $\nu(A) = 0$, under the assumptions of Theorem 1 we also have

$$\nu\left(\left\{\omega \in \Omega : \liminf_n \frac{\sum_{j=1}^n X_j(\omega)}{n} < \int X_1\,d\nu\right\}\right) = 0$$

and

$$\nu\left(\left\{\omega \in \Omega : \limsup_n \frac{\sum_{j=1}^n X_j(\omega)}{n} > -\int -X_1\,d\nu\right\}\right) = 0.$$

In other words, with zero probability empirical averages will eventually lie outside the interval $[\int X_1\,d\nu, -\int -X_1\,d\nu]$.

Second, when $\nu$ is additive we have $\int X_1\,d\nu = -\int -X_1\,d\nu$, and so in this case our result reduces to a standard Kolmogorov limit law

$$\nu\left(\left\{\omega \in \Omega : \lim_n \frac{\sum_{j=1}^n X_j(\omega)}{n} = \int X_1\,d\nu\right\}\right) = 1.$$

On the other hand, when $\nu$ is not additive in some cases it may happen (see [13]) that

$$\nu\left(\left\{\omega \in \Omega : \liminf_n \frac{\sum_{j=1}^n X_j(\omega)}{n} < \limsup_n \frac{\sum_{j=1}^n X_j(\omega)}{n}\right\}\right) = 1.$$

Finally, as anticipated, the closest existing theorem is due to [13]. Our result is more general since [13] assumes that $\Omega$ is compact, $\nu$ is continuous, the r.v.s $X_n$ are continuous, independent and that they satisfy some further technical conditions. Moreover, the proof we provide is different and much simpler. In fact, here we develop a technique that relies on the relations between totally monotone capacities and correspondences, thus making it possible to use existing laws of large numbers for correspondences. This approach might be useful in establishing further generalizations of limit laws to the framework of capacities. This will be the object of future research, along with the possibility of weakening some of the continuity conditions assumed in Theorem 1.

**4. Proof and related material.** Denote by $\mathcal{K}_\Omega$ (resp. $\mathcal{G}_\Omega$) the class of all nonempty compact subsets (resp. open subsets) of $\Omega$; for the sake of completeness, write $\mathcal{B}_\Omega$ instead of $\mathcal{B}$. If $d$ is a Polish metric on $\Omega$, then $\mathcal{K}_\Omega$ is a Polish space when endowed with the Hausdorff metric

$$d_\mathcal{H}(K,L) \equiv \max\left(\max_{k\in K}\min_{l\in L} d(k,l), \max_{l\in L}\min_{k\in K} d(l,k)\right).$$

The Borel $\sigma$-algebra on $\mathcal{K}_\Omega$ is also generated by the class $\{K \in \mathcal{K}_\Omega : K \subseteq G\}_{G\in\mathcal{G}_\Omega}$.

4.1. *Measurable correspondences and totally monotone capacities.* Let $(I,\mathcal{C},\lambda)$ be a nonatomic and complete probability space. A (*compact valued*) *correspondence* $F: I \rightrightarrows \Omega$ is a map with domain $I$ and whose values are nonempty compact subsets of $\Omega$. For any $A \subseteq \Omega$, we put

$$F_{-1}(A) \equiv \{s \in I : F(s) \subseteq A\}.$$

A correspondence $F: I \rightrightarrows \Omega$ is *measurable* if $F_{-1}(G) \in \mathcal{C}$ for every $G \in \mathcal{G}_\Omega$. As well known (see, e.g., [12]), the following facts are equivalent:

- $F$ is measurable;
- $F_{-1}(B) \in \mathcal{C}$ for every $B \in \mathcal{B}_\Omega$;
- $F$ is measurable as a function $F: I \to \mathcal{K}_\Omega$.

When a measurable correspondence $F$ is regarded as a measurable function $F: I \to \mathcal{K}_\Omega$, we denote by $F^{-1}$ its standard inverse image, that is,

$$F^{-1}(\mathcal{E}) \equiv \{s \in I : F(s) \in \mathcal{E}\} \qquad \forall \mathcal{E} \subseteq \mathcal{K}_\Omega,$$

and by $\sigma(F)$ the $\sigma$-algebra generated by $F$, that is,

$$\sigma(F) \equiv \{F^{-1}(\mathcal{D}) : \mathcal{D} \in \mathcal{B}_{\mathcal{K}_\Omega}\}.$$

In the sequel we will need the next lemma, whose standard proof is omitted.

LEMMA 2. *Let $F:I \rightrightarrows \Omega$ be a measurable correspondence. Then*

$$\sigma(F) = \sigma(\{F_{-1}(G) : G \in \mathcal{G}_\Omega\})$$

*and $\{F_{-1}(G) : G \in \mathcal{G}_\Omega\}$ is a $\pi$-class containing $I$.*

A measurable function $f:I \to \Omega$ induces a *probability distribution $P_f$* on $\mathcal{B}_\Omega$ defined by

$$P_f(B) \equiv \lambda(f^{-1}(B)) \qquad \forall B \in \mathcal{B}_\Omega.$$

In a similar way, a measurable correspondence $F:I \rightrightarrows \Omega$ induces a *lower distribution $\nu_F$* on $\mathcal{B}_\Omega$ defined by

$$\nu_F(B) \equiv \lambda(F_{-1}(B)) \qquad \forall B \in \mathcal{B}_\Omega.$$

The next result, which links totally monotone capacities and lower distributions, is essentially due to Choquet [4] (see also [15, 16] and [3]).

LEMMA 3. *A set function $\nu: \mathcal{B}_\Omega \to [0,1]$ is a totally monotone capacity if and only if there exists a measurable correspondence $F:I \rightrightarrows \Omega$ such that $\nu = \nu_F$.*

A *measurable selection* of a correspondence $F:I \rightrightarrows \Omega$ is a measurable function $f:I \to \Omega$ such that $f(s) \in F(s)$ for almost all $s \in I$. The set of all measurable selections of $F$ is denoted by $\operatorname{Sel} F$. The *Aumann integral* (see [2]) of a correspondence $F:I \rightrightarrows \mathbb{R}$ with respect to $\lambda$ is defined by

$$\int F \, d\lambda \equiv \left\{ \int f \, d\lambda : f \in \operatorname{Sel} F \text{ and } f \text{ integrable} \right\}.$$

If $X:\Omega \to \mathbb{R}$ is continuous or simple, and $F$ is a correspondence, then $(X \circ F)(s) \equiv X(F(s))$ is a correspondence [i.e., $X(F(s)) \in \mathcal{K}_\mathbb{R}$ for all $s \in I$]. Moreover, since

$$(X \circ F)_{-1}(A) = F_{-1}(X^{-1}(A)) \qquad \forall A \subseteq \mathbb{R},$$

$X \circ F$ is measurable provided $X$ and $F$ are measurable.

LEMMA 4. *Let $F:I \rightrightarrows \Omega$ be a measurable correspondence, and $X:\Omega \to \mathbb{R}$ be either bounded and continuous or simple and measurable. Then,*

$$\int (X \circ F) \, d\lambda = \left[ \int X \, d\nu_F, - \int -X \, d\nu_F \right].$$

This is an immediate consequence of [3], Theorem 4.1.

4.2. *Proof of Theorem* 1. Suppose first that (ii) holds. By Lemma 3, there exists a measurable correspondence $F: I \rightrightarrows \Omega$ such that $\nu = \nu_F$.

Next we show that the measurable correspondences $\{X_n \circ F\}_{n \geq 1}$ are pairwise independent and identically distributed when regarded as measurable functions $X_n \circ F : I \to \mathcal{K}_\mathbb{R}$.

Let $n, m \geq 1$ and $G_n, G_m \in \mathcal{G}_\mathbb{R}$, then

$$\begin{aligned}
\lambda((X_n \circ F)_{-1}(G_n) &\cap (X_m \circ F)_{-1}(G_m)) \\
&= \lambda(F_{-1}(X_n^{-1}(G_n)) \cap F_{-1}(X_m^{-1}(G_m))) \\
&= \lambda(F_{-1}(X_n^{-1}(G_n) \cap X_m^{-1}(G_m))) \\
&= \nu(X_n^{-1}(G_n) \cap X_m^{-1}(G_m)) \\
&= \nu(X_n^{-1}(G_n))\nu(X_m^{-1}(G_m)) \\
&= \lambda(F_{-1}(X_n^{-1}(G_n)))\lambda(F_{-1}(X_m^{-1}(G_m))) \\
&= \lambda((X_n \circ F)_{-1}(G_n))\lambda((X_m \circ F)_{-1}(G_m)).
\end{aligned}$$

This proves pairwise independence, since for all $j = n, m$, $\{(X_j \circ F)_{-1}(G)\}_{G \in \mathcal{G}_\mathbb{R}}$ is a $\pi$-class containing $I$ and generating the $\sigma$-algebra $\sigma(X_j \circ F)$ (see Lemma 2).

Moreover, for each $n, m \geq 1$, and each open subset $G \in \mathcal{G}_\mathbb{R}$,

$$\begin{aligned}
\lambda((X_n \circ F)^{-1}(\{K \in \mathcal{K}_\mathbb{R} : K \subseteq G\})) \\
&= \lambda(\{X_n \circ F \in \{K \in \mathcal{K}_\mathbb{R} : K \subseteq G\}\}) \\
&= \lambda((X_n \circ F)_{-1}(G)) \\
&= \lambda(F_{-1}(X_n^{-1}(G))) \\
&= \nu(X_n^{-1}(G)) \\
&= \nu(X_m^{-1}(G)) \\
&= \lambda((X_m \circ F)^{-1}(\{K \in \mathcal{K}_\mathbb{R} : K \subseteq G\})).
\end{aligned}$$

This proves identical distribution since $\{K \in \mathcal{K}_\mathbb{R} : K \subseteq G\}_{G \in \mathcal{G}_\mathbb{R}}$ is a $\pi$-class containing $\mathcal{K}_\mathbb{R}$ and generating $\mathcal{B}_{\mathcal{K}_\mathbb{R}}$.

Clearly, for each $n \geq 1$ and each $h \in \operatorname{Sel} X_n \circ F$, $\int h \, d\lambda$ is finite ($h$ is bounded); moreover, by Lemma 4, $\int X_n \circ F \, d\lambda \in \mathcal{K}_\mathbb{R}$.

In sum, $\{X_n \circ F\}_{n \geq 1}$ are pairwise independent and identically distributed measurable correspondences with $\int X_n \circ F \, d\lambda \in \mathcal{K}_\mathbb{R}$ for all $n \geq 1$. A generalization due to [10] (see also [5] and [17]) of a result of Artstein and Vitale [1] guarantees that

$$\lambda\left(\left\{s \in I : \frac{1}{n}\sum_{j=1}^n X_j(F(s)) \to \int X_1 \circ F \, d\lambda\right\}\right) = 1.$$

By Lemma 4,
$$\int X_1 \circ F \, d\lambda = \left[\int X_1 \, d\nu, -\int -X_1 \, d\nu\right].$$

Let $a_n(\omega) = \frac{\sum_{j=1}^n X_j(\omega)}{n}$ and set

$$S_1 \equiv \left\{s \in I : \frac{1}{n}\sum_{j=1}^n X_j(F(s)) \to \left[\int X_1 \, d\nu, -\int -X_1 \, d\nu\right]\right\},$$

$$S_2 \equiv \left\{s \in I : \int X_1 \, d\nu \leq \liminf_n a_n(\omega)\right.$$
$$\left. \leq \limsup_n a_n(\omega) \leq -\int -X_1 \, d\nu \ \forall \omega \in F(s)\right\},$$

$$\Omega_2 \equiv \left\{\omega \in \Omega : \int X_1 \, d\nu \leq \liminf_n a_n(\omega) \leq \limsup_n a_n(\omega) \leq -\int -X_1 \, d\nu\right\}.$$

We want to show that $\nu(\Omega_2) = 1$. Notice that
$$\nu(\Omega_2) = \lambda(\{s \in I : F(s) \subseteq \Omega_2\}) = \lambda(S_2).$$

The next claim will be used to show that $S_1 \subseteq S_2$.

CLAIM 1. *Let $\{K_n\}$ be a sequence in $\mathcal{K}_\mathbb{R}$ such that $K_n \to [\alpha, \beta]$. Then,*
$$\alpha \leq \liminf_n k_n \leq \limsup_n k_n \leq \beta$$
*for each sequence $\{k_n\}$ in $\mathbb{R}$ such that $k_n \in K_n$ for all $n \geq 1$.*

PROOF. By definition of Hausdorff metric, $K_n$ converges to $[\alpha, \beta]$ if and only if $\max(\max_{t_n \in K_n} \min_{r \in [\alpha,\beta]} |t_n - r|, \max_{r \in [\alpha,\beta]} \min_{t_n \in K_n} |r - t_n|) \to 0$, in particular

(1) $$\max_{t_n \in K_n} \min_{r \in [\alpha,\beta]} |t_n - r| \to 0.$$

Let $\{k_{n_j}\}$ be a subsequence of $\{k_n\}$ such that $k_{n_j} \to \ell \in [-\infty, +\infty]$. If $\ell \notin [\alpha, \beta]$, then there exists $\varepsilon > 0$ such that eventually $|k_{n_j} - r| > \varepsilon$ for all $r \in [\alpha, \beta]$. Hence, we have eventually $\min_{r \in [\alpha,\beta]} |k_{n_j} - r| > \varepsilon$, thus contradicting (1). □

If $s \in S_1$, then $\frac{1}{n}\sum_{j=1}^n X_j(F(s)) \to [\int X_1 \, d\nu, -\int -X_1 \, d\nu]$. Hence, for all $\omega \in F(s)$, we have $a_n(\omega) = \frac{1}{n}\sum_{j=1}^n X_j(\omega) \in \frac{1}{n}\sum_{j=1}^n X_j(F(s))$; by Claim 1,
$$\int X_1 \, d\nu \leq \liminf_n a_n(\omega) \leq \limsup_n a_n(\omega) \leq -\int -X_1 \, d\nu.$$

Therefore, $S_1 \subseteq S_2$ and so $\nu(\Omega_2) = \lambda(S_2) \geq \lambda(S_1) = 1$. This completes the proof of the result when (ii) holds.

As to (i), denoting by $\tau$ the Polish topology on $\Omega$, there exists a Polish topology $\tau^* \supseteq \tau$ on $\Omega$ such that $\sigma(\tau^*) = \mathcal{B}_\Omega$, and such that all the $X_n$s are $\tau^*$-continuous (see, e.g., [21]). Since $\nu$ is continuous, then it is a totally monotone capacity with respect to the topology $\tau^*$; we can thus assume that (ii) holds.

**Acknowledgments.** We thank Steve Lalley (the Editor) and an anonymous referee for helpful suggestions.

ISTITUTO DI METODI QUANTITATIVI
 AND IGIER
UNIVERSITÀ BOCCONI
VIALE ISONZO 25
20135 MILANO
ITALY
E-MAIL: fabio.maccheroni@unibocconi.it
URL: http://web.econ.unito.it/gma/fabio.htm

DIPARTIMENTO DI STATISTICA E
 MATEMATICA APPLICATA AND ICER
UNIVERSITÀ DI TORINO
PIAZZA ARBARELLO 8
10122 TORINO
ITALY
E-MAIL: massimo.marinacci@unito.it
URL: http://web.econ.unito.it/gma/massimo.htm